\documentclass[journal,twoside,web]{IEEEtran}
\usepackage{cite}
\usepackage{amsmath,amssymb,amsfonts}
\usepackage{algorithmic}
\usepackage{graphicx}
\usepackage{dsfont}
\usepackage{pstricks}
\usepackage{afterpage}
\usepackage{comment}
\usepackage{cite}
\allowdisplaybreaks[1]
\usepackage{textcomp}

\begin{document}
\title{Prescribed-Time Regulation of Nonlinear Uncertain Systems with Unknown Input Gain and Appended Dynamics}
\author{Prashanth Krishnamurthy and Farshad Khorrami
  \thanks{The authors are with the Control/Robotics Research Laboratory (CRRL), Dept. of ECE, NYU Tandon School of Engineering, Brooklyn, NY 11201, USA (e-mails: prashanth.krishnamurthy@nyu.edu, khorrami@nyu.edu).}
  }

\maketitle

\begin{abstract}
  The prescribed-time stabilization problem for a general class of nonlinear systems with unknown input gain and appended dynamics (with unmeasured state) is addressed. Unlike the asymptotic stabilization problem, the prescribed-time stabilization objective requires convergence of the state to the origin in a finite time that can be arbitrarily picked (i.e., prescribed) by the control system designer irrespective of the initial condition of the system. The class of systems considered is allowed to have general nonlinear uncertain terms throughout the system dynamics as well as uncertain appended dynamics (that effectively generate a time-varying non-vanishing disturbance signal input into the nominal system). The control design is based on a time scale transformation, dynamic high-gain scaling, and adaptation dynamics with temporal forcing terms.
\end{abstract}

\bstctlcite{IEEEexample:BSTcontrol}

\section{Introduction}
While the stabilization/regulation objective typically considered in control designs \cite{KKK95,JK95,Isi99,Kha01} is formulated in term of asymptotic convergence (as time $t\rightarrow\infty$) of the state to a desired state value (e.g., the origin), the control objective of ``finite-time'' stabilization \cite{finite_time1,finite_time2,finite_time3,finite_time4,finite_time5,finite_time6,finite_time7,finite_time8,finite_time9,finite_time10,APA08}
addresses the possibility of achieving the desired convergence properties over a finite time interval. The length of this finite time interval that is attained depends, in general, on the system dynamics and the initial conditions.
Requiring this finite time interval to be a constant that is independent of the initial condition, i.e., requiring that the convergence should be attained within a fixed {\em terminal time} that is independent of initial condition, yields the stronger 
control objective of ``fixed-time'' stabilization \cite{Pol12,finite_time8,PEP16,ZPEP18b,ZPEP18,AGJSD19}.
Further requiring that the fixed finite time should be a parameter that can be arbitrarily ``prescribed'' by the control designer irrespective of the initial condition yields the even stronger control objective of  ``prescribed-time'' stabilization \cite{prescribed_time1,prescribed_time2,JSGL17,JSL17,TYS18,BVAD18,SGLL18,prescribed_time3,SDM19,ZSW19,KKK19a,KKK19b,KKK19c,KK20a,KK20b}.

To force convergence within the specified finite prescribed time, two general prescribed-time stabilizing controller design approaches that have been addressed in the literature can be viewed as state scaling or time scaling:
\begin{itemize}
\item State scaling by a time-dependent function \cite{prescribed_time1,prescribed_time2,prescribed_time3}: By, for example, scaling the state $x$ to define $\tilde x = \mu(t)x$ where $\mu(t)$ is a ``blow-up'' function defined such that $\mu(t)\rightarrow\infty$ as $t\rightarrow T$, a control design that keeps $\tilde x$ bounded will implicitly make $x$ go to 0 as $t\rightarrow T$.
\item Time scaling using a nonlinear temporal transformation \cite{KKK19a,KKK19b,KKK19c,KK20a,KK20b}: Define, for example, $\tau = a(t)$ with $a$ being a function defined such that $a(0)=0$ and $\lim_{t\rightarrow T}a(t)=\infty$. Since this time scale transformation maps $t\in[0,T)$ to $\tau\in[0,\infty)$, a control design that achieves asymptotic convergence in terms of the time variable $\tau$ implicitly achieves prescribed-time convergence in terms of the time variable $t$.
\end{itemize}
The state scaling approach has been applied in \cite{prescribed_time1,prescribed_time2,prescribed_time3} to design prescribed-time stabilizing controllers for classes of systems such as chains of integrators with uncertainties matched with the control input (i.e., normal form).
Prescribed-time stabilizing controllers have been designed for nonlinear strict-feedback-like systems \cite{KKK19a,KKK19b,KKK19c,KK20a,KK20b} using the time scaling approach to convert the prescribed-time stabilization problem into an asymptotic stabilization problem in terms of the transformed time variable and applying the dual dynamic high gain scaling based observer-controller design techniques \cite{KK04f,KK04e,KK07h,KK07f,KK07a,KK08,KK13a,KK16b} to achieve asymptotic stabilization in terms of the transformed time variable. While \cite{KKK19a} considered the prescribed-time stabilization problem under state feedback, output feedback was addressed in \cite{KKK19b}. The adaptation of the control design techniques from \cite{KK04f,KK04e,KK07h,KK07f,KK07a,KK08,KK13a,KK16b} that were originally developed in the context of asymptotic stabilization to the prescribed-time context necessitated introduction of time-dependent forcing terms into the high-gain scaling parameter dynamics and a set of modifications in the controller design and the Lyapunov analysis to achieve prescribed-time convergence instead of asymptotic convergence. 
Uncertain nonlinear systems with general structures of uncertain functions throughout the system dynamics including combinations of unknown parameters (without requiring any known magnitude bounds on unknown parameters) and unmeasured state variables were addressed in \cite{KKK19c} and a dynamic output-feedback prescribed-time stabilizing controller was developed.
A partial state-feedback prescribed-time stabilizing controller was designed for systems with uncertainties in the input gain and non-vanishing input-matched disturbances in addition to uncertain terms throughout the system dynamics in \cite{KK20a}. An output-feedback prescribed-time stabilizing controller was designed for systems with time delays of unknown magnitude in \cite{KK20b}.  

Based on the prescribed-time stabilizing control design in \cite{KK20a}, we consider in this paper a general class of nonlinear systems that include an unknown input gain and time-varying non-vanishing disturbances generated by an uncertain appended dynamics in addition to nonlinear time-varying uncertain terms throughout the system dynamics. The uncertain terms in the system dynamics are allowed to contain both parametric and functional uncertainties without requiring magnitude bounds on the uncertain parameters.  Specifically, we consider a class of nonlinear systems of the following form\footnote{Throughout, a dot above a symbol denotes the derivative with respect to the time $t$, e.g., $\dot x_i=\frac{dx_i}{dt}$. The derivative with respect to the transformed time variable $\tau$ that will be introduced as part of the control design will be written explicitly as, for example, $\frac{dx_i}{d\tau}$.}:
\begin{align}
\dot x_i&=\phi_i(z,x,u,t)+\phi_{(i,i+1)}(x,t)x_{i+1} \ \ , i=1,\ldots,n-1 \nonumber\\
\dot x_n&=\phi_n(z,x,u,t) +h(z,x,u,t)u \nonumber\\
\dot z &= q(z,x,u,t)
\label{system}
\end{align}
where  $x=[x_1,\ldots,x_n]^T\in{\mathcal R}^n$ is the state of the nominal system, $z=[z_1,\ldots,z_{n_z}]^T\in {\mathcal R}^{n_z}$ is the state of an appended dynamics coupled with the $x$ subsystem, and $u\in{\mathcal R}$ is the input\footnote{${\mathcal R}$, ${\mathcal R}^+$, and ${\mathcal R}^k$ denote the set of real numbers, the set of non-negative real numbers, and the set of real $k$-dimensional column vectors, respectively.}. 
$\phi_{(i,i+1)}, i=1,\ldots,n-1$ are known scalar real-valued
continuous functions.  
$\phi_i, i=1,\ldots, n$, $h$, and $q$ are time-varying
scalar real-valued uncertain
functions of their arguments.
The state $x$ of the nominal system is measured while the state $z$ of the appended dynamics is assumed to be unmeasured.
The uncertain function $h$ represents the unknown control input gain, which is allowed to be time-varying and state-dependent.
Furthermore, while $h$ is assumed to have known sign (without loss of generality, assumed positive) and lower-bounded in magnitude by a non-zero constant (to ensure controllability), the lower bound is not required to be known  unlike \cite{KK20a}.
The bounds imposed on functions $\phi_i(z,x,u,t)$ in the assumptions on the system structure in Section~\ref{sec:formulation} allow these functions to depend nonlinearly on the entire system state $x$ as well as appearance of an uncertain parameter $\theta$ (without requirement for a known magnitude bound) and coupling with the state of the appended dynamics. Furthermore, the bound on $\phi_n(z,x,u,t)$ is allowed to contain an additive term that is not required to go to 0 when the state approaches the origin, i.e., a non-vanishing disturbance. A known upper bound on this non-vanishing disturbance is not required unlike \cite{KK20a}. It will be seen in the control design and stability and convergence analysis in Sections~\ref{sec:design} and \ref{sec:proof} that prescribed-time stabilization can be attained for the considered class of systems through several novel ingredients including non-smooth components in the control law, temporal forcing terms in the 
 adaptation dynamics and the scaling parameter dynamics, interconnections between the adaptation dynamics and the scaling parameter dynamics taking into account the time scale transformation $t\rightarrow \tau$, and analysis of the closed-loop system properties.

This paper is organized as follows. The assumptions imposed on the system \eqref{system} are provided in Section~\ref{sec:formulation}. The control design is presented in Section~\ref{sec:design}. The main result of the paper is presented in Section~\ref{sec:proof}.
 Concluding remarks are contained in Section~\ref{sec:conclusion}.

\section{Notations, Control Objective, and Assumptions}
\label{sec:formulation}

\noindent{\bf Notations:}
\begin{itemize}
\item The notation $|a|$ denote Euclidean norm of a vector $a$ or absolute value of a scalar $a$.
The notation $||M||$ denotes Frobenius norm of a matrix $M$.
\item
  The notation $\mbox{diag}(T_1,\ldots,T_m)$ denotes an $m\times m$ diagonal
  matrix with diagonal elements $T_1,\ldots,T_m$.
  Also, $\mbox{lowerdiag}(T_1,\ldots,T_{m-1})$ and $\mbox{upperdiag}(T_1,\ldots,T_{m-1})$ denote the $m\times m$ matrices
  with the lower diagonal entries (i.e., entries at locations $(i+1,i)$ for $i=1,\ldots,m-1$) and upper diagonal entries (i.e., entries at locations $(i,i+1)$ for $i=1,\ldots,m-1$), respectively, being $T_1,\ldots,T_{m-1}$ and zeros
  everywhere else.
  \item $I_m$ denotes the $m\times m$
  identity matrix.
\item The maximum and minimum eigenvalues of a symmetric positive-definite matrix $P$ are denoted by
  $\lambda_{max}(P)$ and $\lambda_{min}(P)$, respectively.
\item The notations $\max(a_1,\ldots,a_n)$ and $\min(a_1,\ldots,a_n)$ indicate the largest and smallest values, respectively, among numbers $a_1,\ldots,a_n$.
\item
Given a vector $a=[a_1,\ldots,a_m]^T$, the notation $|a|_e$ denotes the vector comprised of element-wise magnitudes of the elements of $a$, i.e., $|a|_e=[|a_1|,\ldots,|a_m|]^T$.
Given two vectors $a=[a_1,\ldots,a_m]^T$ and $b=[b_1,\ldots,b_m]^T$, the relation $a \leq_e b$ indicates the set of element-wise inequalities between the corresponding elements of $a$ and $b$, i.e., $|a_i|\leq_e |b_i|, i=1,\ldots,m$.
\item
  Given a scalar $\delta$, the notation $S(\delta)$ denotes the sign of $\delta$, i.e., $S(\delta) = 1$ if $\delta \geq 0$ and $S(\delta) = -1$ otherwise.
\end{itemize}

With $T > 0$ being any prescribed constant,
the control objective is to design a dynamic control law for $u$ using measurement of the signal $x$ so that $x(t)\rightarrow 0$ as $t\rightarrow T$ while also ensuring that $z(t)$ and $u(t)$ remain uniformly bounded over the time interval $t\in[0,T)$, i.e., $\sup_{t\in[0,T)}|z(t)| < \infty$ and $\sup_{t\in[0,T)}u(t) < \infty$.

\noindent{\bf Assumption A1} ({\em lower boundedness away from zero of ``upper diagonal'' terms $\phi_{(i,i+1)}$}):
The inequalities
 $|\phi_{(i,i+1)}(x,t)|\geq \sigma ,1\leq i\leq n-1$ are satisfied for all $x\in{\mathcal R}^n$ and $t\geq 0$ with $\sigma$ being a positive constant. Since $\phi_{(i,i+1)}$ are continuous functions, this assumption can, without loss of generality,
 be stated as
 $\phi_{(i,i+1)}(x,t)\geq \sigma>0 ,1\leq i\leq n-1$.

 \noindent{\bf Assumption A2} ({\em Bounds on uncertain functions $\phi_i$}):
The functions $\phi_i, i=1,\ldots,n$, can be bounded as
\begin{align}
|\phi_i(z,x,u,t)|
  &\leq\Gamma(x_1)\theta
    \sum_{j=1}^i \phi_{(i,j)}(x,t)|x_j|
    \nonumber\\
  &\quad \mbox{ for } i=1,\ldots,n-1
\\
  |\phi_n(z,x,u,t)| &\leq \Gamma(x_1)
                      \Big\{\theta\sum_{j=1}^n \phi_{(i,j)}(x,t)|x_j| + \phi_{n0}\Big\}
\end{align}
for all $x\in{\mathcal R}^n$, $z\in{\mathcal R}^{n_z}$, $u\in{\mathcal R}$, and $t\geq 0$ 
 where $\Gamma(x_1)$ and
$\phi_{(i,j)}(x,t)$ for $i=1,\ldots,n, j=1,\ldots,i$, are known continuous
non-negative functions, and
$\phi_{n0}$ and $\theta$ are unknown non-negative constants.
Positive constants $\epsilon_{(i,j)},i=1,\ldots,n,j=1,\ldots,i$, and
$\tilde\epsilon_{(i,2)},i=2,\ldots,n$, are known
such that $\forall x\in{\mathcal R}^n$ and $t\geq 0$,
\begin{align}
&  \frac{\phi_{(i,1)}(x,t)}{\phi_{(1,2)}(x,t)} \leq \epsilon_{(i,1)} ,\, i=1,\ldots,n
                \nonumber\\
  &  \frac{\phi_{(i,2)}(x,t)}{\sqrt{\phi_{(1,2)}(x,t)\phi_{(2,3)}(x,t)}}\leq \tilde\epsilon_{(i,2)} ,\,
    i=2,\ldots,n
    \nonumber\\
&  \frac{\phi_{(i,j)}(x,t)}{\phi_{(2,3)}(x,t)}\leq \epsilon_{(i,j)}  \ ,
                     i=2,\ldots,n-1, j=2,\ldots,i.
                      \label{a2_2}
\end{align}

\noindent{\bf Assumption A3} ({\em Bound on uncertain input gain $h$}):
The uncertain function $h$ is lower bounded in magnitude by a positive constant $\underline h$ that is not required to be known. Since $h$ is a continuous function, this assumption can, without loss of generality, be stated as  
$h(z,x,u,t)\geq \underline h > 0$ 
for all $x\in{\mathcal R}^n$, $z\in{\mathcal R}^{n_z}$, $u\in{\mathcal R}$, and $t\geq 0$.

\noindent{\bf Assumption A4} ({\em Cascading dominance of ``upper diagonal'' terms $\phi_{(i,i+1)}, i=2,\ldots,n$,}):
Positive constants $\overline\rho_i$ exist such that 
$\phi_{(i,i+1)}(x,t)\geq \overline\rho_i \phi_{(i-1,i)}(x,t)\,,\,\,\,
i=3,\ldots,n-1$, $\forall x\in
{\mathcal R}^n$ and $t\geq 0$.

\noindent{\bf Assumption A5} ({\em Cascading dominance between ``upper diagonal'' terms $\phi_{(1,2)}$ and $\phi_{(2,3)}$}):
Continuous non-negative functions $\overline\phi_{(1,2)}(x_1)$ and $\tilde\phi_{(1,2)}(x_1)$ exist
such that
\begin{align}
  \tilde\phi_{(1,2)}(x_1) \leq \frac{\phi_{(1,2)}(x,t)}{\phi_{(2,3)}(x,t)} \leq \overline\phi_{(1,2)}(x_1)
\end{align}
for all $x\in {\mathcal R}^n$ and $t\geq 0$.

\noindent{\bf Assumption A6} ({\em Assumption on appended dynamics $z$}):
The appended system with the state $z$ and the input $(x,u,t)$ is a bounded-input-bounded-state (BIBS) stable system.

\noindent{\bf Remark 1:}
The Assumptions A1, A4, and A5 are similar to \cite{KKK19a}. The Assumption~A2 is more general than the structure of the corresponding assumption in \cite{KKK19a}. While \cite{KKK19a} assumed bounds of the form $|\phi_i| \leq \Gamma(x_1) \sum_{j=1}^i \phi_{(i,j)}(x)|x_j|$, Assumption~A2 above allows an uncertain parameter $\theta$ (for which no magnitude bounds are required to be known), an additional term $\Gamma(x_1)\phi_{n0}$ in the bound on $|\phi_n|$, and time dependence of $\phi_{(i,j)}$. The additional $\Gamma(x_1)\phi_{n0}$ term in the bound on $|\phi_n|$ allows for the possibility of uncertain non-vanishing disturbance inputs that are driven by the appended dynamics with unmeasured state $z$ and are also time-varying. The Assumptions~A3 and A6 do not have corresponding analogous assumptions in \cite{KKK19a}. Assumption~A3 relates to the unknown input gain $h$ that appears multiplied with the control input $u$ in the system dynamics \eqref{system}. While \cite{KKK19a} assumed that the input $u$ appears with a known gain as $\mu_0(x)u$ with a known function $\mu_0$, the class of systems considered here are allowed to contain an uncertain time-varying state-dependent input gain $h(z,x,u,t)$. Assumption~A3 on this unknown input gain requires only a known lower bound $\underline h$ on $h$ and does not require an upper bound. Assumption~A6 relates to the appended dynamics $z$, which were not considered in \cite{KKK19a}. The role of the appended dynamics here is as a forcing function coupled with various uncertain terms in the system dynamics including $\phi_i$ and $h$.
Also, the Assumptions A2 and A3 are weaker compared to the earlier conference version \cite{KK20a} of this paper. While \cite{KK20a} required the upper bound $\phi_{n0}$ on the non-vanishing part of the uncertain function $\phi_n$ and the lower bound $\underline h$ on the input gain $h$ to be known constants, this requirement is relaxed in this paper. Removing the need for these constants to be known requires several modifications in the control design; specifically, while the design of $u_1$ in \eqref{sfcontrol} utilized $\underline h$ in \cite{KK20a} and the design of $u_2$ in \eqref{u2defn} utilized $\underline h$ and $\phi_{n0}$, additional time-dependent functions are introduced in this paper in place of $\underline h$ and $\phi_{n0}$. By designing these time-dependent forcing functions appropriately in combination with various modifications in the stability analysis, the need for knowledge of constants $\phi_{n0}$ and $\underline h$ is removed. 
$\diamond$

\section{Control Design}
\label{sec:design}

\subsection{Design of Control Law $u$}
The control input $u$ is designed as comprised of two components:
\begin{align}
  u &= u_1 + u_2
      \label{u_combined}
\end{align}
where $u_1$ defined below is picked based on a pair of coupled Lyapunov inequalities as discussed in Section~\ref{sec:coupled_lyap} and $u_2$ is designed as part of Section~\ref{sec:design_freedoms} based on a Lyapunov analysis taking into account the various uncertain terms in the system dynamics including the unknown input gain $h$.
The first component $u_1$ in \eqref{u_combined} is designed as
\begin{align}
  u_1&=-\frac{r^n}{\gamma_1(t)}K_c \eta
     \label{sfcontrol}
\end{align}
where:
\begin{itemize}
\item
  $r$ is a dynamic high-gain scaling parameter whose dynamics to be designed in Section~\ref{sec:design_freedoms} will be such that $r$ is a monotonically non-decreasing signal in time. $r$ will be initialized such that $r(0)\geq 1$; hence, $r(t)\geq 1$ for all time $t$.
\item $\gamma_1:{\mathcal R}^+\rightarrow {\mathcal R}^+$ is a function that will be designed as part of Section~\ref{sec:design_freedoms}.
\item $\eta=[\eta_2,\ldots,\eta_n]^T$ with $\eta_i$ being scaled state variables defined as:
  \begin{align}
    \eta_2 &= \frac{x_2 + \zeta(x_1,\hat\theta)}{r}
             \ ; \ \eta_i = \frac{x_i}{r^{i-1}}, i=3,\ldots,n.
             \label{eta2defn_etaidefn}
  \end{align}
  \item
    $\zeta
    $ is a function defined to be of the form
  \begin{align}
    \zeta(x_1,\hat\theta) &= \hat\theta x_1\zeta_1(x_1)
                   \label{zetadefn}
  \end{align}
  with $\zeta_1$ being a function that will be designed as part of Section~\ref{sec:design_freedoms}
  and $\hat\theta$ is a dynamic adaptation parameter. The dynamics of $\hat\theta$ will also be designed as part of Section~\ref{sec:design_freedoms} and will be such that $\hat\theta$ is a monotonically non-decreasing signal as a function of time. $\hat\theta$ will be initialized such that $\hat\theta(0) \geq 1$; hence, $\hat\theta(t)\geq 1$ for all time $t$.
\item
  $K_c=[k_2,\ldots,k_n]$ with $k_i,i=2,\ldots,n$, being functions of $(x,t)$ that will be designed below.
\end{itemize}
The dynamics of the scaled state vector $\eta$ defined above under the control law given by
\eqref{u_combined} and \eqref{sfcontrol} can be written as\footnote{For notational convenience, we drop the arguments of
  functions whenever no confusion will result.}
\begin{align}
  \dot\eta&=rA_c \eta -\frac{\dot r}{r}D_c\eta+\Phi+H\eta_2+\Xi
            \nonumber\\
  &\quad
            - B\frac{h - \gamma_1(t)}{\gamma_1(t)}rK_c\eta + B h\frac{u_2}{r^{n-1}}
\label{sfetamatdyn}
\end{align}
where  $A_c$  is the matrix of dimension $(n-1)\times (n-1)$ in which the $(i,j)^{th}$
element is given by
$A_{c_{(i,i+1)}}=\phi_{(i+1,i+2)}, i=1,\ldots, n-2$,
$A_{c_{(n-1,j)}}=-k_{j+1},  j=1,\ldots,n-1$,
and zeros elsewhere.
$B$ is the column vector of length $(n-1)$ of form $[0,\ldots,0,1]^T$.
Also,
\begin{align}
D_c&=\mbox{diag}(1,2,\ldots,n-1)
\ ; \ 
\Phi=[\frac{\phi_2}{r},\ldots,\frac{\phi_n}{r^{n-1}}]^T
      \label{Phidefn}
      \\
H&=[\hat\theta[\zeta_1'(x_1)x_1+\zeta_1]\phi_{(1,2)},0,\ldots,0]^T
\\
  \Xi&=[\frac{(\phi_1-\zeta\phi_{(1,2)})\hat\theta[\zeta_1'(x_1)x_1+\zeta_1] + \dot{\hat\theta}x_1\zeta_1(x_1)}{r},
\nonumber\\
       &\quad 0,\ldots,0]^T.
\end{align}
where
$\zeta_1'(x_1)$ denotes the partial derivative 
of the function $\zeta_1$ evaluated at $x_1$.

If the input gain $h$ were a known function, $h$ could be used in place of $\gamma_1(t)$ in \eqref{sfcontrol} to cancel out the function $h$ in the resulting dynamics of $\eta$, i.e., to remove the term involving $(h-\gamma_1(t))$ in dynamics \eqref{sfetamatdyn}. However, since $h$ is an uncertain function and even the lower bound $\underline h$ is unknown, the function $\gamma_1(t)$ is introduced in \eqref{sfcontrol} and will be designed in Section~\ref{sec:design_freedoms} to handle  the ``mismatch'' term involving $(h-\gamma_1(t))$.

\subsection{Time Scale Transformation}
\label{sec:temporal_scaling}

Define a time scale transformation $\tau=a(t)$ satisfying the following conditions:
\begin{itemize}
\item
  $a$ is a twice continuously differentiable monotonically increasing function over $[0,T)$ with
$a(0)=0$ and $\lim_{t\rightarrow T}a(t)=\infty$.
\item
  Denoting $a'(t)=\frac{da}{dt}$, a positive constant $a_0$ exists such that $a'(t)\geq a_0$ for all $t\in[0,T)$. The first condition and this condition imply that the function $a$ is invertible.
  Denote the inverse function by $a^{-1}$, i.e., $t = a^{-1}(\tau)$.
\item
  Denoting the function $a'(t)=\frac{da}{dt}$ expressed in terms of the new time variable $\tau$ by $\alpha(\tau)$, i.e., $\alpha(\tau)=a'(a^{-1}(\tau))$,
the function $\alpha(\tau)$ grows at most polynomially as $\tau\rightarrow\infty$, i.e., a polynomial $\overline\alpha(\tau)$ exists such that $\alpha(\tau)\leq \overline{\alpha}(\tau)$ for all $\tau\in[0,\infty)$. Also, $\frac{d\alpha}{d\tau}$ grows at most polynomially as $\tau\rightarrow\infty$.
\end{itemize}

\noindent{\bf Remark 2:}
The time scale transformation $\tau = a(t)$ maps the finite time interval $[0,T)$ in terms of the original time variable $t$ to the infinite time interval $[0,\infty)$ in terms of the transformed time variable $\tau$. Hence, the prescribed-time control objective formulated as convergence objectives as $t\rightarrow T$ are equivalent to the analogous convergence objectives as $\tau\rightarrow\infty$. From the definition of the time scale transformation, we have
\begin{align}
  dt &= \frac{1}{\alpha(\tau)} d\tau.
       \label{dtau_dt}
\end{align}
As noted in \cite{KKK19a,KKK19b}, an infinite number of functions exist that satisfy the conditions required on the function $a$ above. For example, one choice for the function $a$ is
$a(t) = \frac{a_0t}{1-\frac{t}{T}}$ with $a_0$ being any positive constant. With this choice of the function $a$, we have
$a'(t) = \frac{a_0}{(1-\frac{t}{T})^2}$ and
$\alpha(\tau) = a_0 \bigg(\frac{\tau}{a_0 T} + 1\bigg)^2$. Note that
$\alpha(\tau)$ and $\frac{d\alpha}{d\tau}$ do indeed grow at most polynomially with the time $\tau$ as required in the conditions introduced above on $\alpha$.
In the stability and convergence analysis in Section~\ref{sec:proof}, it will be seen that this polynomial growth condition is indeed crucial in showing that the high-gain scaling parameter $r$ grows at most polynomially with the time $\tau$; since $x_i$ can be written in terms of combinations of $\eta_i$ and powers of $r$, it will be seen that the polynomial growth property of $r$ is crucial in inferring that exponential convergence of $\eta$ to 0 as $\tau\rightarrow\infty$ implies exponential convergence of $x_i$ to 0.  Similarly, since the control law for $u$ involves terms comprising of combinations such as $r^n\eta$ as seen in \eqref{sfcontrol}, we will see in the analysis in Section~\ref{sec:proof} that the polynomial growth property of $r$ is crucial in inferring convergence to 0 of these terms in the control law from the exponential convergence of $\eta$ to 0. 
$\diamond$

\subsection{Lyapunov Function}
Define
\begin{align}
V&=\frac{1}{2}x_1^2+r\eta^T P_c\eta
\label{Vdefn}
\end{align}
where $P_c$ is a constant symmetric positive-definite matrix that will be defined in Section~\ref{sec:coupled_lyap} based on the solution of a pair of coupled Lyapunov inequalities. Differentiating \eqref{Vdefn}
and using the property that $dt = \frac{d\tau}{\alpha(\tau)}$, we have
\begin{align}
 \frac{dV}{d\tau}
  &=\frac{1}{\alpha(\tau)}\Bigg\{x_1[\phi_1+(r\eta_2 -\zeta)\phi_{(1,2)}]
    \nonumber\\
  &\quad
    +r^2\eta^T[P_c A_c +A_c^T P_c]\eta
  +2r\eta^T P_c(\Phi+H\eta_2+\Xi)
    \nonumber\\
  &\quad
    -2r^2\frac{h-\gamma_1(t)}{\gamma_1(t)}\eta^T P_c B K_c \eta
  + 2\eta^T P_c B h\frac{u_2}{r^{n-2}}
  \Bigg\}
    \nonumber\\&\quad
  -\frac{dr}{d\tau} \eta^T[P_c \tilde D_c+ \tilde D_c P_c]\eta
  \label{Vdot1}
\end{align}
where $\tilde D_c=D_c- \frac{1}{2}I_{n-1}$.

\subsection{Coupled Lyapunov Inequalities}
\label{sec:coupled_lyap}
Assumption A3 is the {\em cascading dominance} condition (among the upper diagonal terms $\phi_{(2,3)},\ldots,\phi_{(n-1,n)}$) introduced in \cite{KK04f}; under this condition, it was shown in \cite{KK04f,KK06} that a constant symmetric positive-definite matrix $P_c$ and a function $K_c(x,t) = [k_2(x,t),\ldots,k_n(x,t)]$ (whose elements appear in the definition of the matrix $A_c$) can be constructed  such that the following coupled Lyapunov inequalities are satisfied with some positive constants $\nu_{c}$, $\underline\nu_{c}$,
and $\overline\nu_{c}$:
\begin{align}
  P_c A_c+A_c^T P_c \leq -\nu_c \phi_{(2,3)} I
  \nonumber\\
  \underline \nu_c I \leq P_c \tilde D_c + \tilde D_c P_c\leq \overline\nu_c I.
  \label{coupled_lyap}
\end{align}

\subsection{Inequality Bounds on Terms Appearing in Lyapunov Inequality \eqref{Vdot1}}
\label{sec:ineq_bounds}

Using the bounds on uncertain terms $\phi_i$ in Assumption A2, the definition of $\Phi$ in \eqref{Phidefn}, the definitions of the scaled state variables $\eta_i$ in \eqref{eta2defn_etaidefn}, and the property that $r\geq 1$, we have
\begin{align}
  \Phi &\leq_e \theta\Gamma \Phi_1 \frac{|x_1|}{r} + \theta\Gamma \Phi_M |\eta|_e + \theta\Gamma \Phi_2 \frac{\hat\theta|\zeta_1x_1|}{r} + \Gamma B \frac{\phi_{n0}}{r^{n-1}}
\end{align}
where
\begin{itemize}
\item $\Phi_1 = [\phi_{(2,1)},\ldots,\phi_{(n,1)}]^T$; $\Phi_2 = [\phi_{(2,2)},\ldots,\phi_{(n,2)}]^T$
\item $\Phi_M$ is the matrix of dimension $(n-1)\times (n-1)$ with $(i,j)^{th}$ element $\phi_{(i+1,j+1)}$ for $i=1,\ldots,n-1, j\leq i$, and zeros everywhere else
\end{itemize}
Hence (with some conservative overbounding for algebraic simplicity),
\begin{align}
  2r\eta^T P_c \Phi &\leq
  2\lambda_{max}(P_c)|\eta|\theta\Gamma\bigg\{
                     |\Phi_1| |x_1| + |\Phi_2| \hat\theta|\zeta_1 x_1|
                     \nonumber\\
  &\quad
                     + r |\Phi_M| |\eta| 
  \bigg\} + 2\frac{\theta\Gamma}{r^{n-2}}|\eta^TP_c|_e B \phi_{n0}
  \\
                   &\leq
                     \frac{2\theta^2\Gamma^2}{\zeta_0\phi_{(1,2)}}\lambda_{max}^2(P_c)|\eta|^2\{|\Phi_1|^2 + |\Phi_2|^2 \hat\theta^2 \zeta_1^2\}
                     \nonumber\\
  &\quad
                     + \zeta_0\phi_{(1,2)}x_1^2
    + 2r\theta\Gamma\lambda_{max}(P_c)|\Phi_M||\eta|^2
    \nonumber\\
                   &\quad
                     + 2\frac{\Gamma}{r^{n-2}}|\eta^TP_c|_e B \phi_{n0}
                     \label{bnd_ineq0}
\end{align}
with $\zeta_0$ being any positive constant and with $|\eta^TP_c|_e$ denoting the vector comprised of the element-wise magnitudes of the elements of the vector $\eta^T P_c$ as per the notation defined in Section~\ref{sec:formulation}.

Using Assumption A2, the other terms in the Lyapunov inequality \eqref{Vdot1} can also be upper bounded as:
\begin{align}
  x_1\phi_1 &\leq \theta \Gamma x_1^2  \phi_{(1,1)}
              \label{bnd_ineq1}
  \\
  x_1r\eta_2\phi_{(1,2)} &\leq \frac{\nu_c}{4}r^2 \phi_{(2,3)} |\eta|^2
+\frac{1}{\nu_c}x_1^2\frac{\phi_{(1,2)}^2}{\phi_{(2,3)}}
  \\
  2r\eta^T P_cH\eta_2 &\leq 2\hat\theta r \lambda_{max}(P_c)\phi_{(1,2)}|\zeta_1'x_1+\zeta_1||\eta|^2
\\
  2r\eta^T P_c\Xi &\leq \zeta_0 \phi_{(1,2)} x_1^2
                 \! +\! \frac{2}{\zeta_0\phi_{(1,2)}}\lambda_{max}^2(P_c) |\eta|^2 \Big[\!\dot{\hat\theta}^2\zeta_1^2
\nonumber\\
            &\quad              + (\theta\Gamma\phi_{(1,1)}+|\zeta_1|\hat\theta\phi_{(1,2)})^2
\nonumber\\
&\quad             \times (\zeta_1'x_1+\zeta_1)^2\hat\theta^2 \Big]
              \label{bnd_ineqn}
\end{align}
Using the inequalities in \eqref{coupled_lyap} and \eqref{bnd_ineq0}--\eqref{bnd_ineqn}, \eqref{Vdot1} yields
\begin{align}
  \frac{dV}{d\tau}
  &\leq \frac{1}{\alpha(\tau)}\Bigg\{
    -x_1^2\hat\theta\zeta_1\phi_{(1,2)}
    -\frac{3}{4}\nu_c \phi_{(2,3)} r^2|\eta|^2
    \nonumber\\
  &\quad
    + q_1(x_1)\phi_{(1,2)}x_1^2
    + \theta^* q_2(x_1)\phi_{(1,2)} x_1^2
    \nonumber\\
  &\quad
    +rw_1\Big(x_1,\hat\theta,\frac{\dot{\hat\theta}}{\phi_{(1,2)}}\Big)\phi_{(1,2)}|\eta|^2
    \nonumber\\
  &\quad
    +r\theta^* w_2(x_1,\hat\theta)\phi_{(2,3)}|\eta|^2
    -2r^2\frac{h-\gamma_1(t)}{\gamma_1(t)}\eta^T P_c B K_c \eta
    \nonumber\\
  &\quad
    + 2\eta^T P_c B h\frac{u_2}{r^{n-2}}
    + 2\frac{\Gamma}{r^{n-2}}|\eta^TP_c|_e B\phi_{n0} 
    \Bigg\}
    \nonumber\\
    &\quad
    -\underline\nu_c\frac{dr}{d\tau} |\eta|^2
  \label{Vdot3}
\end{align}
where $\theta^* = (1+\theta+\theta^2)$ is an uncertain positive constant and
\begin{align}
  q_1(x_1) &= \frac{1}{\nu_c}\overline\phi_{(1,2)}(x_1)
           + 2\zeta_0 
           \label{q1defn}
  \\
  q_2(x_1) &= \Gamma(x_1)\epsilon_{(1,1)}
             \label{q2defn}
             \\
  w_1\Big(x_1,\hat\theta,\frac{\dot{\hat\theta}}{\phi_{(1,2)}}\Big) &=
                                         2\hat\theta \lambda_{max}(P_c)|\zeta_1'(x_1)x_1+\zeta_1(x_1)|
                                         \nonumber\\
  &\quad 
                                         +
                                         \frac{2}{\zeta_0}\lambda_{max}^2(P_c)\Big(\frac{\dot{\hat\theta}}{\phi_{(1,2)}(x_1)}\Big)^2\zeta_1^2(x_1)
    \nonumber\\
           &\quad 
                                         +
             \frac{4}{\zeta_0}\lambda_{max}^2(P_c)\zeta_1^2(x_1)
             \nonumber\\
  &\quad\times
             \hat\theta^4(\zeta_1'(x_1)x_1+\zeta_1(x_1))^2
                                         \label{w1defn}
  \\
  w_2(x_1,\hat\theta) &=
                        2\lambda_{max}^2(P_c)\frac{\Gamma^2(x_1)}{\zeta_0}
                        \nonumber\\
           &\quad
             \times
                        \Bigg[\overline\phi_{(1,2)}(x_1)\sqrt{\sum_{i=2}^n\epsilon_{(i,1)}^2}
                                         \nonumber\\
           &\quad
                                         + \hat\theta^2\zeta_1^2(x_1)\sqrt{\sum_{i=2}^n\tilde\epsilon_{(i,2)}^2} \Bigg]
\nonumber\\
           &\quad
             + 2 \Gamma(x_1)\lambda_{max}(P_c)\sqrt{\sum_{i=2}^n\sum_{j=2}^i\epsilon_{(i,j)}^2}
             \nonumber\\
           &\quad
             + \frac{4}{\zeta_0}\epsilon_{(1,1)}^2\overline\phi_{(1,2)}(x_1)\Gamma^2(x_1)\lambda_{max}^2(P_c)
             \nonumber\\
           &\quad\ 
             \times (\zeta_1'(x_1)x_1+\zeta_1(x_1))^2\hat\theta^2
           \label{w2defn}
\end{align}
Note that the functions $q_1(x_1)$, $q_2(x_1)$ $w_1\Big(x_1,\hat\theta,\frac{\dot{\hat\theta}}{\phi_{(1,2)}(x,t)}\Big)$, and $w_2(x,\hat\theta,\dot{\hat\theta})$ involve only known functions and quantities.
The third argument of the definition of $w_1$ is written in terms of the combination $\frac{\dot{\hat\theta}}{\phi_{(1,2)}(x,t)}$  rather than as simply $\dot{\hat\theta}$ separately since
it will be seen (in Lemma 2 in Section~\ref{sec:proof}) that it can be shown that $\frac{\dot{\hat\theta}}{\phi_{(1,2)}(x,t)}$ grows at most polynomially as a function of the time $\tau$ and that this property can then be used to show (in Lemma~3 in Section~\ref{sec:proof}) that $r$ grows at most polynomially as a function of the time $\tau$.

\subsection{Designs of Functions $\zeta_1$ and $\gamma_1$, Dynamics of $r$ and $\hat\theta$, and Control Law Component $u_2$}
\label{sec:design_freedoms}
The design freedoms appearing in the right hand side of \eqref{Vdot3} are $\zeta_1$, $\gamma_1$, $\frac{dr}{d\tau}$, and $u_2$. In addition, the dynamics of $\hat\theta$, i.e., $\frac{d\hat\theta}{d\tau}$ is also a design freedom as will be seen below.
The function $\zeta_1$ is designed such that the negative $x_1^2\hat\theta\zeta_1\phi_{(1,2)}$ term in the right hand side of \eqref{Vdot3} dominates over the positive $q_1\phi_{(1,2)}x_1^2$ and $\theta^*q_2\phi_{(1,2)}x_1^2$ terms, but with the unknown constant $\theta^*$ replaced by $\hat\theta$, which is a dynamic adaptation state variable. Hence, noting that $\hat\theta \geq 1$, we pick the function $\zeta_1$ such that
\begin{align}
  \frac{1}{4}\zeta_1(x_1) &= \max\bigg\{ \underline\zeta, q_1(x_1)+q_2(x_1) \bigg\}
\label{zeta1defn}
\end{align}
with any constant $\underline\zeta>0$.

To design the dynamics of the high-gain scaling parameter $r$, we use the basic motivation from
the dynamic high-gain scaling control designs for asymptotic stabilization (e.g., \cite{KK04f}) that the dynamics of $r$ should be designed such that the time derivative of $r$ is ``large enough'' (in a nonlinear function sense) until $r$ itself becomes ``large enough'' (also in a nonlinear function sense). Furthermore, the state-dependent form of these two ``large enough'' functions should be designed based on Lyapunov analysis such that desired Lyapunov inequalities hold both under the case that the time derivative of $r$ is large enough and the case that $r$ is large enough. For this purpose, the 
dynamics of $r$ are designed to be of the form
\begin{align}
  \frac{dr}{d\tau}&= \lambda\bigg(R\Big(x_1,\hat\theta,\frac{\dot{\hat\theta}}{\phi_{(1,2)}(x,t)}\Big) + \alpha(\tau)-r\bigg)
  \nonumber\\
  &\quad
  \times
  [\Omega(r,x,\hat\theta,\dot{\hat\theta},t) + \tilde\alpha(\tau)]
\nonumber\\
&\qquad \mbox{ with } r(0)\geq \max\{1, \alpha(0)\}
                    \label{rdot}
\end{align}
where $\tilde\alpha(\tau)$ denotes $\frac{d\alpha}{d\tau}$.
Here, $\lambda:{\mathcal R}\rightarrow{\mathcal R}^+$ can be picked to be any
continuous function such that
$\lambda(s)=1$ for $s\geq 0$ and $\lambda(s)=0$ for
$s\leq-\epsilon_r$ where $\epsilon_r$ can be picked to be any positive constant. With such a choice of the function $\lambda$, it is seen that $\frac{dr}{d\tau}$ is ``large'' (i.e., $\frac{dr}{d\tau} = \Omega +\tilde\alpha$) when $r$ is relatively small and on the other hand, when $r$ becomes ``large'' (i.e., $r\geq R + \alpha + \epsilon_r$), $\frac{dr}{d\tau}$ goes to 0.
The functions $R$ and $\Omega$ are picked as
\begin{align}
  R\bigg(x_1,\hat\theta,\frac{\dot{\hat\theta}}{\phi_{(1,2)}(x,t)}\bigg) &= \max\bigg\{1,
  \nonumber\\
  &\!\! 
  \frac{4}{\nu_c}\Big[w_1\Big(x_1,\hat\theta,\frac{\dot{\hat\theta}}{\phi_{(1,2)}(x,t)}\Big)\overline\phi_{(1,2)}(x_1)
  \nonumber\\
  &\!\! 
  +\hat\theta w_2(x_1,\hat\theta)\Big]
  \bigg\}
             \label{Rdefn}
\end{align}
\begin{align}
  \Omega(r,x,\hat\theta,\dot{\hat\theta},t) &= \frac{r}{\underline\nu_c a_0}\bigg[w_1(x_1,\hat\theta,\dot{\hat\theta})\phi_{(1,2)}(x,t)
                  \nonumber\\
           &\quad                           + \hat\theta w_2(x,t)\phi_{(2,3)}(x,t)\bigg].
             \label{Omegadefn}
\end{align}
The function $R$ is chosen such that when $r\geq R$, the negative term involving $\nu_c\phi_{(2,3)}r^2|\eta|^2$ in the right hand side of \eqref{Vdot3} dominates over the positive $rw_1\phi_{(1,2)}|\eta|^2$ and $r\theta^*w_2\phi_{(2,3)}|\eta|^2$ terms, but with the unknown constant $\theta^*$ replaced by the adaptation parameter $\hat\theta$. The function $\Omega$ is chosen such that when $\frac{dr}{d\tau} \geq\Omega$, the negative term involving $\underline\nu_c \frac{dr}{d\tau}|\eta|^2$ in the right hand side of \eqref{Vdot3} dominates over the positive $rw_1\phi_{(1,2)}|\eta|^2$ and $r\theta^*w_2\phi_{(2,3)}|\eta|^2$ terms, but again with the unknown constant $\theta^*$ replaced by $\hat\theta$. Hence, effectively, when $r$ is relatively small (i.e., when $r < R$), the derivative $\frac{dr}{d\tau}$ is relatively large (i.e., $\frac{dr}{d\tau}\geq \Omega$) by the form of the dynamics of $r$ in \eqref{rdot} and therefore the the negative term involving $\underline\nu_c \frac{dr}{d\tau}|\eta|^2$ in the right hand side of \eqref{Vdot3} dominates over the positive $rw_1\phi_{(1,2)}|\eta|^2$ and $r\hat\theta w_2\phi_{(2,3)}|\eta|^2$ terms. On the other hand, when $r$ is sufficiently large (i.e., when $r\geq R$) the negative term involving $\nu_c\phi_{(2,3)}r^2|\eta|^2$ in the right hand side of \eqref{Vdot3} dominates over the positive $rw_1\phi_{(1,2)}|\eta|^2$ and $r\hat\theta w_2\phi_{(2,3)}|\eta|^2$ terms.

As noted in Section~\ref{sec:ineq_bounds}, the appearance of $\dot{\hat\theta}$ in the dynamics of $r$ is written in terms of the combination $\frac{\dot{\hat\theta}}{\phi_{(1,2)}}$ rather than simply $\dot{\hat\theta}$ since this combination can be shown (Lemma~2 in Section~\ref{sec:proof}) to grow at most polynomially as a function of time $\tau$, a property that will be used in showing (Lemma~3 in Section~\ref{sec:proof}) that $R$ and therefore $r$ grow at most polynomially as a function of $\tau$.

With the dynamics of $r$ as designed in \eqref{rdot}--\eqref{Omegadefn}, it is seen that $\frac{dr}{d\tau} \geq 0$ for all time $\tau \geq 0$.
Also, noting that $r$ is initialized such that $r(0)\geq\alpha(0)$ and noting from \eqref{rdot} that we have $\frac{dr}{d\tau} \geq \tilde\alpha(\tau) = \frac{d\alpha}{d\tau}$ at any time instant at which $r \leq R + \alpha(\tau)$ where $R\geq 0$ from \eqref{Rdefn}, we see that $r \geq \alpha(\tau)$ for all time $\tau$ in the maximal interval of existence of solutions.

The dynamics of the adaptation parameter $\hat\theta$ is designed as (the motivation for this form of the dynamics of $\hat\theta$ can be seen from the augmented Lyapunov function $\overline V$ in \eqref{Vbardefn} and its derivative \eqref{Vbardot}):
\begin{align}
  \frac{d\hat\theta}{d\tau} &= \tilde\alpha(\tau) + c_\theta
  \frac{\chi(r,x,\hat\theta,t)}{\alpha(\tau)}
               \label{thetahatdyn}
                     \ \mbox{ with } \hat\theta(0) \geq \max\{1, \alpha(0)\}
\end{align}
where 
\begin{align}
  \chi(r,x,\hat\theta,t) &=
                           \phi_{(1,2)}(x,t)q_2(x_1)x_1^2
                           \nonumber\\
  &\quad
                     + r w_2(x_1,\hat\theta)\phi_{(2,3)}(x,t)|\eta|^2 
               \label{chidefn}
\end{align}
where $c_\theta$ is any positive constant.
From \eqref{thetahatdyn}, it is seen that $\hat\theta(0)\geq \alpha(0)$ 
and also $\frac{d\hat\theta}{d\tau} \geq \frac{d\alpha(\tau)}{d\tau}$ for all times $\tau$. 
Hence,  $\hat\theta \geq \alpha(\tau)$ at all times $\tau$ in the maximal interval of existence of solutions.

Considering the remaining terms in the right hand side of \eqref{Vdot3}, i.e., the terms involving $r^2\frac{h-\gamma_1(t)}{\gamma_1(t)}\eta^T P_c B K_c \eta$ and $\frac{\Gamma}{r^{n-2}}|\eta^T P_c|_e B\phi_{n0}$, the control law component $u_2$ is designed such that the term $\eta^T P_c B h \frac{u_2}{r^{n-2}}$ in the right hand side of \eqref{Vdot3} dominates over these two terms, but with a time-dependent function $\gamma_2(t)$ in place of $\phi_{n0}$ since $\phi_{n0}$ is unknown. The function
$\gamma_2:{\mathcal R}^+\rightarrow {\mathcal R}^+$ will be designed below.
Hence,
the component $u_2$ of the control input signal as defined in \eqref{u_combined} is designed as:
\begin{align}
  u_2 &= -S(\eta^T P_c B) \bigg\{|K_c\eta|r^n \bigg[\frac{1}{\gamma_1(t)} +\hat\theta_1\bigg]
        \nonumber\\
  &\quad
        + \Gamma(x_1)\bigg[\frac{\gamma_2(t)}{\gamma_1(t)}+\hat\theta_1\bigg]\Bigg\}
  \label{u2defn}
\end{align}
where $\hat\theta_1$ is an adaptation parameter whose dynamics will be designed below in \eqref{theta1hatdyn}.
The dynamics of $\hat\theta_1$ will be designed such that $\hat\theta_1$ is a monotonically non-decreasing signal as a function of time and $\hat\theta_1$ will be initialized such that $\hat\theta_1(0)\geq 0$. Hence, $\hat\theta_1(t)\geq 0$ for all time $t$.
In \eqref{u2defn}, the notation $S(\delta)$ with $\delta$ being a scalar denotes the sign of $\delta$ as defined in Section~\ref{sec:formulation}.
Analogous to \eqref{sfcontrol}, the time-dependent function $\gamma_1(t)$ is used in the denominator of multiple terms in \eqref{u2defn} in place of $h$ since the function $h$ is unknown.

Consider the two cases (a) $r\geq R$; (b) $r < R$. Under case (b), we have $\frac{dr}{d\tau} \geq \Omega$ from the form of the dynamics of $r$ in \eqref{rdot} corresponding to the property as discussed above that the dynamics \eqref{rdot} ensures that either $r$ or its derivative $\frac{dr}{d\tau}$ is ``large''. Using \eqref{zeta1defn}--\eqref{u2defn}, it is seen that in both  cases (a) and (b), \eqref{Vdot3} reduces to
\begin{align}
  \frac{dV}{d\tau}
  &\leq \frac{1}{\alpha(\tau)}\Bigg\{
    -\frac{3}{4}x_1^2\hat\theta\zeta_1\phi_{(1,2)}
    -\frac{1}{2}\nu_c \phi_{(2,3)} r^2|\eta|^2
    \Bigg\}
    \nonumber\\
    &\quad 
      + (\theta^* - \hat\theta)\frac{\chi(r,x,\hat\theta,t)}{\alpha(\tau)}
      \nonumber\\
  &\quad
    -\frac{1}{\alpha(\tau)}\Bigg\{
    2r^2|\eta^T P_c B K_c \eta|\bigg[\frac{h}{\gamma_1}+h\hat\theta_1-\frac{|h-\gamma_1(t)|}{\gamma_1(t)}\bigg]
    \nonumber\\
  &\quad
    + 2|\eta^T P_c B|\frac{\Gamma}{r^{n-2}}  \bigg[\frac{h\gamma_2}{\gamma_1}+h\hat\theta_1-\phi_{n0}\bigg]
    \Bigg\}.
  \label{Vdot4}
\end{align}
It was noted above that the dynamics \eqref{thetahatdyn} and \eqref{rdot} for $\hat\theta$ and $r$, respectively, imply that  $r \geq \alpha(\tau)$ and $\hat\theta \geq \alpha(\tau)$ for all time $\tau$. Hence, \eqref{Vdot4} yields
\begin{align}
  \frac{dV}{d\tau}
  &\leq 
    -\frac{3}{4}x_1^2\zeta_1\phi_{(1,2)}
    -\frac{1}{2}\nu_c \phi_{(2,3)} r|\eta|^2
    \nonumber\\
  &\quad 
    + (\theta^* - \hat\theta)\frac{\chi(r,x,\hat\theta,t)}{\alpha(\tau)}
    \nonumber\\
  &\quad
    -\frac{1}{\alpha(\tau)}\Bigg\{
    2r^2|\eta^T P_c B K_c \eta|\bigg[\frac{h}{\gamma_1}+h\hat\theta_1-\frac{|h-\gamma_1(t)|}{\gamma_1(t)}\bigg]
    \nonumber\\
  &\quad
    + 2|\eta^T P_c B|\frac{\Gamma}{r^{n-2}}  \bigg[\frac{h\gamma_2}{\gamma_1}+h\hat\theta_1-\phi_{n0}\bigg]
    \Bigg\}.
    \label{Vdot5}
\end{align}
Therefore, comparing with the definition of $V$ in \eqref{Vdefn}, we have
\begin{align}
  \frac{dV}{d\tau}
  &\leq
    -\kappa V
    + (\theta^* - \hat\theta)\frac{\chi(r,x,\hat\theta,t)}{\alpha(\tau)}
    \nonumber\\
  &\quad
    -\frac{1}{\alpha(\tau)}\Bigg\{
    2r^2|\eta^T P_c B K_c \eta|\bigg[\frac{h}{\gamma_1}+h\hat\theta_1-\frac{|h-\gamma_1(t)|}{\gamma_1(t)}\bigg]
    \nonumber\\
  &\quad
    + 2|\eta^T P_c B|\frac{\Gamma}{r^{n-2}}  \bigg[\frac{h\gamma_2}{\gamma_1}+h\hat\theta_1-\phi_{n0}\bigg]
    \Bigg\}
    \label{Vdot6}
\end{align}
where
\begin{align}
  \kappa &= \min\bigg\{\frac{3\zeta_0\sigma}{2},\frac{\nu_c\sigma}{2\lambda_{max}(P_c)}\bigg\}.
           \label{kappadefn}
\end{align}
Noting that $h$, $\gamma_1$, and $\gamma_2$ are non-negative, noting that $h\geq \underline h$,
and defining $\theta_1*=\max(\frac{1}{\underline h},\frac{\phi_{n0}}{\underline h})$,
\eqref{Vdot6} yields
\begin{align}
  \frac{dV}{d\tau}
  &\leq
    -\kappa V
    + (\theta^* - \hat\theta)\frac{\chi(r,x,\hat\theta,t)}{\alpha(\tau)}
    \nonumber\\
  &\quad
    +
    \underline h
    (\theta_1^*-\hat\theta_1)
    \frac{\chi_1(r,x,\hat\theta,t)}{\alpha(\tau)}
    \label{Vdot7}
\end{align}
where
\begin{align}
  \chi_1(r,x,\hat\theta,t) &= 2r^2|\eta^T P_c B K_c \eta| + 2|\eta^T P_c B|\frac{\Gamma}{r^{n-2}}.
\end{align}
Based on the form of the dynamics in \eqref{Vdot7}, the dynamics of $\hat\theta_1$ are designed as
\begin{align}
  \frac{d\hat\theta_1}{d\tau} &= c_{\theta 1}\frac{\chi_1(r,x,\hat\theta,t)}{\alpha(\tau)}.
                                \label{theta1hatdyn}
\end{align}
The temporal forcing term $\tilde\alpha(\tau)$ is incorporated into the dynamics of $\hat\theta$ in \eqref{thetahatdyn} to ensure that $\hat\theta\geq \alpha(\tau)$, a property that is required to be able to infer \eqref{Vdot5} from \eqref{Vdot4}. 
Noting that $\frac{d}{d\tau}(\hat\theta - \alpha(\tau)) = \chi(r,x,\hat\theta,t)$ from \eqref{thetahatdyn}, it is seen from \eqref{Vdot7} that the signal $(\hat\theta - \alpha(\tau))$ would suffice as the adaptation state variable to address the uncertain parameter $\theta^*$.
Hence, defining an augmented Lyapunov function $\overline V$ that adds to $V$ an additional quadratic component in terms of $(\hat\theta - \alpha(\tau) - \theta^*)$ as well as 
a quadratic component in terms of $(\hat\theta_1-\theta_1^*)$, i.e.,
\begin{align}
  \overline V &= V + \frac{1}{2c_\theta}(\hat\theta - \alpha(\tau) - \theta^*)^2 + \frac{\underline h}{2c_{\theta 1}}(\hat\theta_1 - \theta_1^*)^2,
                \label{Vbardefn}
\end{align}
we have from \eqref{thetahatdyn} and \eqref{Vdot5} and noting that $\chi(r,x,\hat\theta,t)\geq 0$:
\begin{align}
  \frac{d\overline V}{d\tau}
  &\leq 
    -\frac{3}{4}x_1^2\zeta_1\phi_{(1,2)}
    -\frac{1}{2}\nu_c \phi_{(2,3)} r|\eta|^2.
  \label{Vbardot}
\end{align}
While, as we will seen in Section~\ref{sec:proof}, \eqref{Vbardot} can be used to show existence of solutions of the closed-loop dynamical system over the time interval $\tau\in[0,\infty)$, it will not directly enable showing exponential convergence (since the quadratic terms involving the adaptation parameters do not appear on the right hand side of \eqref{Vbardot}). Showing exponential convergence of $x_1$ and $\eta$ to 0 will be crucial in proving closed-loop stability since, for example, the boundedness of $u_1$ will be proved by showing that $r$ grows at most polynomially as a function of time $\tau$ while $\eta$ goes to 0 exponentially. Hence, to show exponential convergence, we will also want to ensure that an inequality of the form $\frac{dV}{d\tau}\leq -\kappa V$ is also satisfied at least after a sub-interval of the overall time interval $\tau\in[0,\infty)$. From \eqref{Vdot6}, we will for this purpose want to ensure that after some finite time, the following inequalities are satisfied:
\begin{align}
  \hat\theta &\geq \theta^*
               \label{thetahatcond}
               \\
  \frac{h}{\gamma_1}+h\hat\theta_1 &\geq \frac{|h-\gamma_1|}{\gamma_1}
                                     \label{theta1hatcond1}
                                     \\
\frac{h\gamma_2}{\gamma_1}+h\hat\theta_1 &\geq \phi_{n0}
                                           \label{theta1hatcond2}
\end{align}
are satisfied. From the dynamics of $\hat\theta$ in \eqref{thetahatdyn}, it will be seen that the inequality
\eqref{thetahatcond}
will be satisfied after some finite time. To ensure that \eqref{theta1hatcond1} and \eqref{theta1hatcond2} are satisfied after some finite time, we pick the functions $\gamma_1$ and $\gamma_2$ such that $\frac{1}{\gamma_1}$ and $\gamma_2$ go to $\infty$ as $t\rightarrow T$, i.e., as $\tau\rightarrow \infty$, by defining
\begin{align}
  \gamma_1(t) &= \frac{1}{c_{\gamma 1}\alpha(a(t)) + \tilde c_{\gamma 1}}
                \label{gamma1defn}
                \\
  \gamma_2(t) &= [c_{\gamma 2}\alpha(a(t))  + \tilde c_{\gamma 2}]\gamma_1(t)
                \label{gamma2defn}
\end{align}
with $c_{\gamma 1}$ and $c_{\gamma 2}$ being any positive constants and $\tilde c_{\gamma 1}$ and $\tilde c_{\gamma 2}$ being any non-negative constants.
From the conditions imposed on the function $\alpha$ and the definitions of the functions $\gamma_1$ and $\gamma_2$ \eqref{gamma1defn} and \eqref{gamma2defn}, it is seen that $\frac{1}{\gamma_1}$ and $\gamma_2$ grow at most polynomially as functions of the time $\tau$.

\section{Stability Analysis and Main Result}
\label{sec:proof}
In this section, a sequence of lemmas is established based on the adaptive controller design in Section~\ref{sec:design}.
Let the maximal interval of existence of solutions of the closed-loop system be $[0,\tau_f)$ in terms of the new time variable $\tau$. From Lemmas~1--4, it is shown that $\tau_f=\infty$, i.e., solutions exist over the infinite time interval $\tau\in[0,\infty)$. Thereafter, various convergence properties are shown in Lemmas 5--7. The main prescribed-time stabilization result of this paper (Theorem~1) is then stated and proved based on the Lemmas 1--7. 

\noindent{\bf Lemma 1:} The signals $V$, $x_1$, $\sqrt{r}|\eta|$, $(\hat\theta - \alpha(\tau))$, and $\hat\theta_1$ are uniformly bounded over $[0,\tau_f)$.

\noindent{\bf Proof of Lemma 1:} From \eqref{Vbardot}, it is seen that $\frac{d\overline V}{d\tau} \leq 0$ implying that $\overline V$ is uniformly bounded over the maximal interval of existence of solutions $[0,\tau_f)$ of the closed-loop system.  From the definitions of $V$ and $\overline V$ in \eqref{Vdefn} and \eqref{Vbardefn}, respectively, the statement of Lemma~1  follows. $\diamond$

\noindent{\bf Lemma 2:} The signals $\hat\theta(a^{-1}(\tau))$ and  $\frac{\dot{\hat\theta}(a^{-1}(\tau))}{\phi_{(1,2)}(x(a^{-1}(\tau)),a^{-1}(\tau))}$ grow at most polynomially in the time variable $\tau = a(t)$ as $\tau \rightarrow \infty$.

\noindent{\bf Proof of Lemma 2:} It was seen as part of Lemma~1 that $(\hat\theta - \alpha(\tau))$ is uniformly bounded over $[0,\tau_f)$. Noting that $\alpha(\tau)$ grows at most polynomially in $\tau$ due to the conditions imposed in Section~\ref{sec:temporal_scaling} on the choice of the function $\alpha$, it follows that $\hat\theta$ grows at most polynomially as a function of time $\tau$. Noting the dynamics of the adaptation variable $\hat\theta$ in \eqref{thetahatdyn} and using the Assumptions A1 and A5, it is seen that
\begin{align}
  \frac{\dot{\hat\theta}}{\phi_{(1,2)}(x,t)} &\leq \alpha(\tau)\tilde\alpha(\tau) + c_\theta
                                               \bigg\{
                                               q_2(x_1)x_1^2
                                               \nonumber\\
                                             &\quad
                                               + \frac{r}{\tilde\phi_{(1,2)}(x_1)} w_2(x_1,\hat\theta)|\eta|^2 
                                               \bigg\}.
                                               \label{thetahatdotbound}
\end{align}
Note that $\alpha(\tau)$ and $\tilde\alpha(\tau)$ grow at most polynomially in $\tau$ by construction (Section~\ref{sec:temporal_scaling}). It was seen in Lemma~1 that $x_1$ and $r|\eta|^2$ are uniformly bounded over $[0,\tau_f)$.
It was noted above that $\hat\theta$ grows at most polynomially as a function of $\tau$. From the definition of $w_2$ in \eqref{w2defn}, $\hat\theta$ appears polynomially (as terms involving $\hat\theta^2$) in $w_2$.
Therefore, it follows from \eqref{thetahatdotbound} that $\frac{\dot{\hat\theta}}{\phi_{(1,2)}(x,t)}$ grows at most polynomially in the time $\tau$. $\diamond$

\noindent{\bf Lemma 3:} The signal $r(a^{-1}(\tau))$ grows at most polynomially in time $\tau$ as $\tau \rightarrow \infty$.

\noindent{\bf Proof of Lemma 3:} Using the Lemmas 1 and 2, it is seen that $w_1\Big(x_1,\hat\theta,\frac{\dot{\hat\theta}}{\phi_{(1,2)}}\Big)$ and $w_2(x_1,\hat\theta)$ defined in \eqref{w1defn} and \eqref{w2defn}, respectively, grow at most polynomially in time $\tau$. Hence, it is seen from \eqref{Rdefn} that $R\Big(x_1,\hat\theta,\frac{\dot{\hat\theta}}{\phi_{(1,2)}(x,t)}\Big)$ grows at most polynomially with time $\tau$.
From \eqref{rdot}, it is seen that $\dot r = 0$ at any time instant $\tau$ at which $r \geq R\Big(x_1,\hat\theta,\frac{\dot{\hat\theta}}{\phi_{(1,2)}(x,t)}\Big) + \alpha(\tau) + \epsilon_r$. 
By the conditions imposed on the function $\alpha(\tau)$ in Section~\ref{sec:temporal_scaling}, $\alpha(\tau)$ is also
polynomially upper bounded in $\tau$.
Hence, 
$R\Big(x_1,\hat\theta,\frac{\dot{\hat\theta}}{\phi_{(1,2)}(x,t)}\Big) + \alpha(\tau) + \epsilon_r$ and therefore
$r$ as well grow at most polynomially as a function of time $\tau$. $\diamond$

\noindent{\bf Lemma 4:} Solutions to the closed-loop dynamical system formed by the given system \eqref{system} and the designed dynamic controller from Section~\ref{sec:design} exist over time interval $\tau\in[0,\infty)$.

\noindent{\bf Proof of Lemma 4:} It is seen from Lemma 1 that $x_1$, $\sqrt{r}\eta$, and $\hat\theta_1$ remain uniformly bounded over $[0,\tau_f)$ while it is seen from Lemmas 2 and 3 that $\hat\theta$ and $r$ grow at most polynomially in $\tau$. Hence, it follows that all closed-loop signals are bounded over any finite time interval $\tau\in[0,\overline\tau_f)$ and therefore solutions to the closed-loop dynamical system exist over the time interval $\tau\in[0,\infty)$, i.e., $\tau_f=\infty
$.  $\diamond$

\noindent{\bf Lemma 5:} A finite constant $\tau_0 \geq 0$ exists such that for all time $\tau \geq \tau_0$, the inequality $\frac{dV}{d\tau}\leq -\kappa V$ is satisfied where the constant $\kappa > 0$ is as defined in \eqref{kappadefn}.

\noindent{\bf Proof of Lemma 5:} From the dynamics of the adaptation parameter $\hat\theta$ in \eqref{thetahatdyn}, it was noted in Section~\ref{sec:design_freedoms} that $\hat\theta\geq \alpha(\tau)$ for all time $\tau$, implying (due to the construction of the function $\alpha(\tau)$ in Section~\ref{sec:temporal_scaling}) that $\hat\theta$ goes to $\infty$ as $\tau\rightarrow\infty$. Hence, a finite constant $\tau_1 > 0$ exists such that
\eqref{thetahatcond} is satisfied
for all time $\tau \geq \tau_1$.
Similarly,
from the construction of $\alpha(\tau)$ and the definitions of $\gamma_1$ and $\gamma_2$, 
finite constants $\tau_2$ and $\tau_3$ exist such that
\eqref{theta1hatcond1}
and 
\eqref{theta1hatcond2}
are satisfied for all times $\tau \geq \tau_2$ and $\tau \geq \tau_3$, respectively.
Hence,
defining $\tau_0=\max(\tau_1,\tau_2,\tau_3)$, it is seen 
from \eqref{Vdot6} that for all times $\tau\geq \tau_0$, we have $\frac{dV}{d\tau}\leq -\kappa V$ with $\kappa$ given in \eqref{kappadefn}. $\diamond$

\noindent{\bf Lemma 6:} The signals $V$, $x_1$, and $\sqrt{r}|\eta|$ go to 0 exponentially as $\tau \rightarrow\infty$. 

\noindent{\bf Proof of Lemma 6:} From Lemma 5, it is seen that a finite constant $\tau_0 > 0$ exists such that for all times $\tau\geq \tau_0$, the inequality $\frac{dV}{d\tau}\leq -\kappa V$ is satisfied. Therefore, $V$ goes to 0 exponentially as $\tau \rightarrow\infty$. From the definition of $V$ in \eqref{Vdefn}, it follows that $x_1$ and $\sqrt{r}|\eta|$ go to 0 exponentially as $\tau\rightarrow\infty$. $\diamond$

\noindent{\bf Lemma 7:} $\eta$ goes to 0 exponentially as $\tau \rightarrow\infty$. Also, $u$ is uniformly bounded over time interval $\tau\in[0,\infty)$.

\noindent{\bf Proof of Lemma 7:} From Lemma 6, we see that $\eta$ goes to 0 exponentially as $\tau\rightarrow\infty$ since $r\geq 1$ for all time $\tau$.
Since, from Lemma 2, $\hat\theta$ grows at most polynomially while from Lemma 6, $x_1$ goes to 0 exponentially, it is seen that $\zeta(x_1,\hat\theta)$ defined in \eqref{zetadefn} goes to 0 exponentially as $\tau\rightarrow\infty$.
Since $r$ grows at most polynomially in time $\tau$ from Lemma 3 while $\eta$ goes to 0 exponentially, it follows from the definition of $\eta_2,\ldots,\eta_n$ in \eqref{eta2defn_etaidefn} that $x_2,\ldots,x_n$ go to 0 exponentially as $\tau\rightarrow\infty$.
Also, $r^n\eta$ goes to 0 exponentially as $\tau\rightarrow\infty$.
Hence, noting that $\frac{1}{\gamma_1}$ grows at most polynomially from \eqref{gamma1defn}, it is seen that
$u_1$ defined in \eqref{sfcontrol} goes to 0 exponentially as $\tau\rightarrow\infty$. Similarly,
noting that $\gamma_2$ also grows at most polynomially while $\hat\theta_1$ remains uniformly bounded, it also follows from the definition of $u_2$ in \eqref{u2defn} that $u_2$ is uniformly bounded over the time interval $\tau\in[0,\infty)$. Therefore, the signal $u=u_1+u_2$ is uniformly bounded over time interval $\tau\in[0,\infty)$. 

\noindent{\bf Theorem 1:} Under the Assumptions A1--A6, the closed-loop dynamical system formed by the given system \eqref{system} and the dynamic controller (of dynamic order 3 -- with state variables $r$, $\hat\theta$, and $\hat\theta_1$) designed in Section~\ref{sec:design} with $T>0$ being arbitrarily picked by the designer satisfies the property that starting from any initial conditions for $x$ and $z$, the signals $x$, $z$, and $u$ satisfy $\lim_{t\rightarrow T}|x(t)| = 0$, $\sup_{t\in [0,T)}|u(t)| < \infty$, and $\sup_{t\in[0,T)}|z(t)| < \infty$.

\noindent{\bf Proof of Theorem 1:}
Noting that $x_2=r\eta_2 - \zeta$ and $x_i=\eta_i r^{i-1}, i=3,\ldots,n$, it follows from the Lemmas 6 and 7 that $x=[x_1,\ldots,x_n]^T$ goes to 0 exponentially as $\tau\rightarrow\infty$.  
From Lemma 7 and Assumption A6, it is seen that $u$ and $z$ are uniformly bounded over time interval $\tau\in[0,\infty)$.
Since $\tau\rightarrow\infty$ corresponds to $t\rightarrow T$, these properties hold as $t\rightarrow T$. 
$\diamond$

\noindent{\bf Remark 3:}
The designed prescribed-time stabilizing adaptive dynamic controller is of dynamic order 3 with the controller state variables being the dynamic scaling parameter $r$ with the dynamics shown in \eqref{rdot}, the adaptation parameter $\hat\theta$ with the dynamics shown in \eqref{thetahatdyn}, and the adaptation parameter $\hat\theta_1$ with the dynamics shown in \eqref{theta1hatdyn}. The overall controller is given by the definition of scaled state vector $\eta$ in \eqref{eta2defn_etaidefn}, control law given by the combination of \eqref{u_combined}, \eqref{sfcontrol}, and \eqref{u2defn},  the choice of the function $\zeta$ in \eqref{zetadefn} and \eqref{zeta1defn},
the choices of the functions $\gamma_1$ and $\gamma_2$ in \eqref{gamma1defn} and \eqref{gamma2defn}, respectively,
the scaling parameter dynamics in \eqref{rdot}, \eqref{Rdefn}, and \eqref{Omegadefn}, and the adaptation parameter dynamics in \eqref{thetahatdyn} and \eqref{theta1hatdyn}.

\noindent{\bf Remark 4:}
As seen in the closed-loop analysis above, several signals in the closed-loop system such as $r$, $\alpha(\tau)$, and $\hat\theta$ go to $\infty$ as $t\rightarrow T$ (with at most polynomial growth as a function of the transformed time variable $\tau$). The polynomial growth of these signals implies that effective control gains go to $\infty$ as $t\rightarrow T$. This is essentially expected since as noted in \cite{prescribed_time1,prescribed_time2,prescribed_time3},  indeed any approach for regulation in finite time (including optimal control designs with a terminal constraint and sliding mode based controllers with time-varying gains) will share the property that effective control gains go to $\infty$ as $t\rightarrow T$. However, it is to be noted that, as proved above, the actual control signal $u$ remains bounded over the time interval $[0,T)$. Also, $x$ goes to 0 as $t\rightarrow T$. Nevertheless, numerical challenges in the implementation of the controller can be posed by the unbounded growth of the effective control gains as $t\rightarrow T$. As noted in \cite{KKK19a,KKK19b,KKK19c},
numerical difficulties can be alleviated using several techniques such as adding a dead zone on the state $x$, adding a saturation on the control gains, implementing the dynamics of the high-gain scaling parameter $r$ via a temporally scaled version $\tilde r = r\tilde\rho(\tau)$, and setting the effective terminal time $\overline T$ in controller implementation to be a constant slightly larger than the desired prescribed time $T$.

\section{Illustrative Example}
\label{sec:example}

Consider the fifth-order system
\begin{align}
  \dot x_1 &= (1+x_1^2)x_2
             \nonumber\\
  \dot x_2 &= (1+x_1^4)x_3 + \theta_a \cos(x_2z_1)x_2
             \nonumber\\
  &\quad
             + \theta_b [1+\cos(tu)]e^{x_1}x_1^2\sin(z_2) 
             \nonumber\\
  \dot x_3 &= \Big[1+\frac{1}{2}\sin(t)\cos(z_2)+x_1^4(1+e^{-|z_1|})\Big]u
             \nonumber\\
           &\quad
             + \theta_cx_1^2 \cos(x_3z_1)x_2 + \theta_d(1+x_1^2)
             \nonumber\\
  \dot z_1 &= -100z_1 + z_2
             \nonumber\\
  \dot z_2 &= -100z_2 + x_3^2 + u
             \label{ex_system}
\end{align}
where $\theta_a$, $\theta_b$, $\theta_c$, and $\theta_d$ are uncertain parameters (with no magnitude bounds required to be known). This system is of the form \eqref{system} with $\phi_{(1,2)}(x_1)=1+x_1^2$, $\phi_{(2,3)}(x_1)=1+x_1^4$, and $h(z,x,u,t)=1+\frac{1}{2}\sin(t)\cos(z_2)+x_1^4(1+e^{-|z_1|})$.
Assumption A1 is satisfied with the constant $\sigma=1$.
 Assumption~A2 is satisfied with
 $\theta=\max\{\theta_a,2\theta_b,\theta_c\}/c_\beta$,
 $\phi_{n0}=\theta_d/c_\beta$,
 $\Gamma(x_1)=c_\beta\max(e^{x_1}|x_1|,1+x_1^2)$,
 $\phi_{(1,1)}=\phi_{(3,1)}=\phi_{(3,3)}=0$, and $\phi_{(2,1)}=\phi_{(2,2)}=\phi_{(3,2)}=1$
 with $c_\beta$ being any positive constant.
 It is seen that inequalities \eqref{a2_2} in Assumption~A2 are trivially satisfied.
 Note that the forms of the various uncertain terms in the dynamics are not required to be known as long as bounds of the form in Assumption~A2 are known to be satisfied.
 Assumption~A3 is satisfied with $\underline h = 0.5$.
 Assumption~A4 is trivially satisfied since $n=3$.
 Assumption~A5 is satisfied with $\overline\phi_{(1,2)}=\frac{3}{2}$ and $\tilde\phi_{(1,2)}=\frac{1+x_1^2}{1+x_1^4}$.
 Noting that the $z$ dynamics is a stable linear system with $x_3$ and $u$ as inputs, it is seen that Assumption~A6 is satisfied.
 Using the constructive procedure in \cite{KK04f,KK04g,KK06} for solution of coupled Lyapunov inequalities,
a symmetric positive-definite matrix $P_c$ and functions $k_2$ and $k_3$ can be found
to satisfy \eqref{coupled_lyap} as $P_c=\tilde a_c\left[\begin{array}{cc}3 & 1 \\ 1&1 \end{array}\right]$, $k_2=5\phi_{(2,3)}$, and $k_3=4\phi_{(2,3)}$, and with $\nu_c=1.675 \tilde a_c$, $\underline\nu_c=\tilde a_c$, and $\overline\nu_c=5\tilde a_c$ with $\tilde a_c$ being any positive constant.
The function $\alpha$ is picked as in Remark~2.
The functions $\gamma_1$ and $\gamma_2$ are picked as in \eqref{gamma1defn} and \eqref{gamma2defn}.
Defining $\eta_2=\frac{x_2+\hat\theta x_1\zeta_1}{r}$ and $\eta_3=\frac{x_3}{r^2}$, we have
  $u_1 = -r^3[k_2\eta_2 + k_3\eta_3]/\gamma_1(t)$ from \eqref{sfcontrol}.
  Since $\phi_{(1,1)}$ and therefore $\epsilon_{(1,1)}$ are 0, we have $q_2=0$ from \eqref{q2defn}. Also, $q_1$ and therefore $\zeta_1$ are constants since $\overline\phi_{(1,2)}$ was found above to be a constant.
  The control component $u_2$ is defined as in \eqref{u2defn} and the overall control input $u$ is defined as $u = u_1+u_2$ from \eqref{u_combined}.
  The dynamics of $r$ are as shown in \eqref{rdot} where the functions $R$ and $\Omega$ are computed following the procedure in Section~\ref{sec:design} and using sharper bounds taking the specific system structure \eqref{ex_system} into account and noting that several terms in the upper bounds vanish since $\phi_{(1,1)}$, etc., are zero for this system and $\zeta_1$ is a constant.
  The dynamics of the adaptation parameters $\hat\theta$ and $\hat\theta_1$ are as shown in \eqref{thetahatdyn} and \eqref{theta1hatdyn}.

\begin{figure*}
  \centerline{\includegraphics[width=0.8\textwidth]{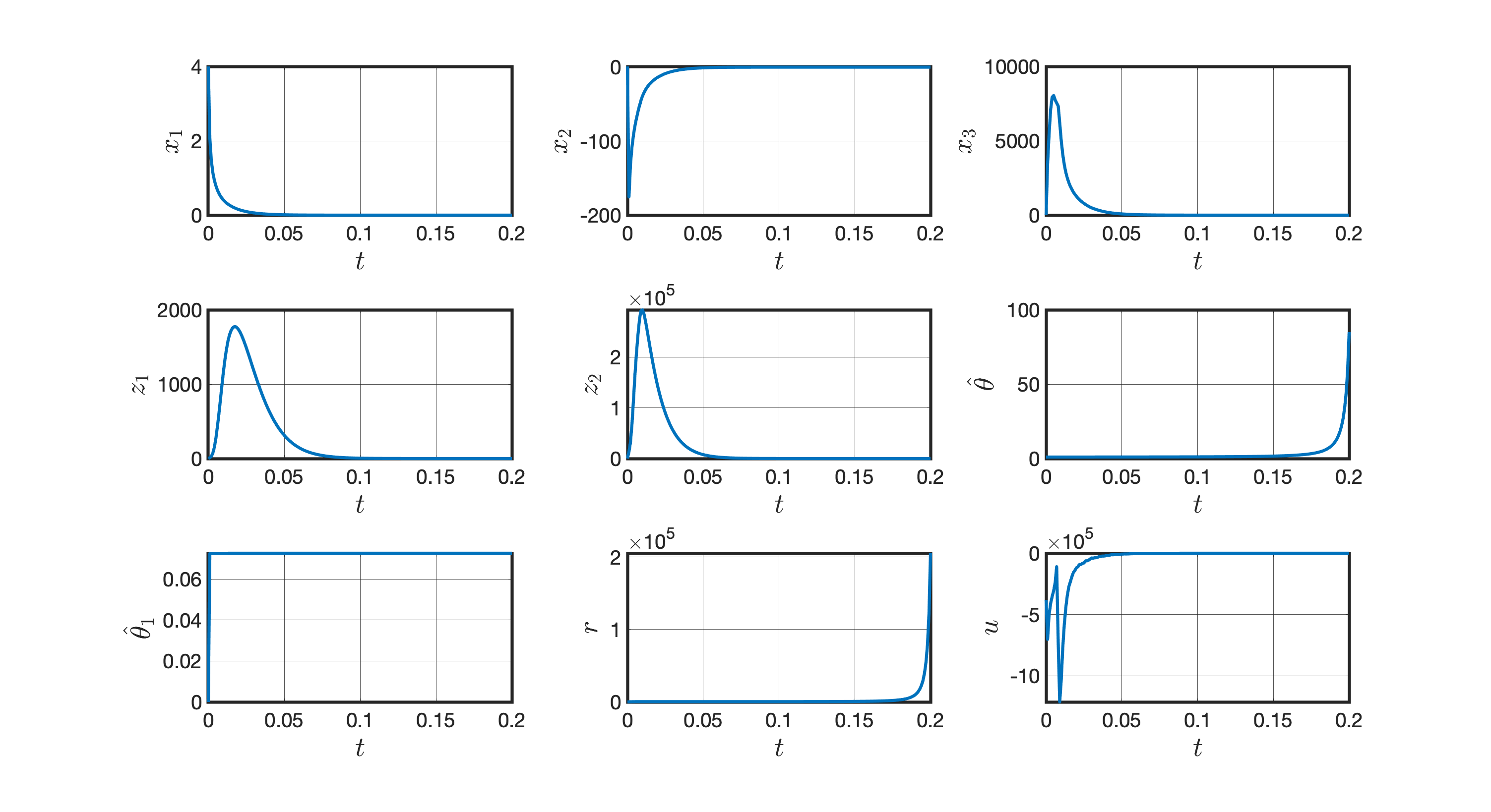}}
  \caption{Simulations for closed-loop system (system \eqref{ex_system} in closed loop with the prescribed-time stabilizing adaptive controller).}
  \label{fig:sim}
\end{figure*}

  The prescribed terminal time is picked to be $T=0.2$~s. To avoid numerical issues as discussed in Remark 4, the effective terminal time $\overline T$ in the controller implementation is defined as $\overline T = 0.205$~s.
  The parameters in the definitions of the time-dependent functions $\alpha$, $\gamma_1$, and $\gamma_2$ are picked as
  $a_0=0.05$, $c_{\gamma 1}=0.01$, $\tilde c_{\gamma 1}=0.5$, and $c_{\gamma 2}=\tilde c_{\gamma 2}=10^{-4}$.
  Also, $\zeta_0=0.25$, $\tilde a_c=0.05$, $c_\theta=10^{-4}$, $c_{\theta 1}=0.01$, and $c_\beta=10^{-4}$.
  The values of the uncertain parameters $\theta_a$, $\theta_b$, $\theta_c$, and $\theta_d$ are picked for simulations as $\theta_a=\theta_b=\theta_c=\theta_d=2$. With the initial conditions for the system state $[x_1,x_2,x_3,z_1,z_2]^T$ specified as $[4,1,1, 1,1]^T$ and the initial conditions for the controller state $[\hat\theta,\hat\theta_1,r]^T$ specified as $[1,0,1]^T$, the closed-loop trajectories and control input signal are shown in Figure~\ref{fig:sim}.

 \section{Conclusion}
 \label{sec:conclusion}

 By combining a non-smooth control component ($u_2$ which involves the sign of $\eta^T P_c B$), time-dependent forcing functions in the definitions of both $u_1$ and $u_2$, an adaptation dynamics that incorporates temporal forcing terms, a time scale transformation $t\rightarrow\tau$, and dynamic scaling-based control design, it was shown that a prescribed-time stabilizing controller can be designed for a general class of nonlinear uncertain systems. The class of nonlinear systems considered allows several types of uncertainties including uncertain input gain and appended dynamics that effectively generate non-vanishing disturbances as well as a general structure of state-dependent uncertain terms throughout the system dynamics. While the adaptation parameter $\hat\theta$ and the dynamic scaling parameter $r$ grow (at most polynomially) as a function of the transformed time variable $\tau$, it was shown that the system state and input remain uniformly bounded and the system state $x$ converges to 0 in the prescribed time irrespective of the initial conditions of the system. Determining if similar control design approaches can be applied to other and more general classes of systems such as general cascade structures and non-triangular and feedforward systems as well as systems with unknown sign of the control gain remain topics for further research.

\bstctlcite{IEEEexample:BSTcontrol}
 \bibliographystyle{IEEEtran}
\bibliography{refs}
\bstctlcite{IEEEexample:BSTcontrol}

\end{document}